\documentclass[10pt,a4paper]{article}
\usepackage{mathrsfs}
\usepackage{epsfig, graphicx}
\usepackage{latexsym,amsfonts,amsbsy,amssymb}
\usepackage{amsmath,amsthm}
\usepackage{color,enumitem}
\usepackage[colorlinks, citecolor=blue]{hyperref}
\usepackage{bm,upgreek}
\usepackage{diagbox}[2011/11/22]
\usepackage{tabularx}
\usepackage{epstopdf}
\makeatletter 

\@addtoreset{equation}{section}
\makeatother 
\allowdisplaybreaks 

\def\w\eta{\widetilde{\eta}}

\newcommand{\bc}{\begin{center}}
\newcommand{\ec}{\end{center}}
\newcommand{\be}{\begin{eqnarray}}
\newcommand{\ee}{\end{eqnarray}}

\newcommand{\ben}{\begin{eqnarray*}}
\newcommand{\een}{\end{eqnarray*}}

\newtheorem{remark}{Remark}[section]
\newtheorem{algorithm}{Algorithm}[section]

\begin{document}
\title{Multilevel correction type of adaptive finite element method for Hartree--Fock equation\footnote{This work was supported by the
Project of Cultivation for young top-motch Talents of Beijing Municipal Institutions (Grant No. BPHR202203022),
the Beijing Municipal Natural Science Foundation (Grant No. 1232001),
and the National Natural Science Foundation of China (Grant No. 12371386). }}
\author{
Fei Xu\footnote{School of Mathematics, Statistics and Mechanics, Beijing University of Technology, Beijing, 100124, P.R. China.
    (fxu@bjut.edu.cn)}  }

\date{} \maketitle

\begin{abstract}
This paper proposes an efficient algorithm for solving the Hartree--Fock equation combining a multilevel correction scheme with an adaptive refinement technique to improve computational efficiency. The algorithm integrates a multilevel correction framework with an optimized implementation strategy.
Within this framework, a series of linearized boundary value problems are solved, and their approximate solutions are corrected by solving small-scale Hartree--Fock equations in low-dimensional correction spaces. The correction space comprises a coarse space and the solution to the linearized boundary value problem, enabling high accuracy while preserving low-dimensional characteristics.
The proposed algorithm efficiently addresses the inherent computational complexity of the Hartree--Fock equation. Innovative correction strategies eliminate the need for direct computation of large-scale nonlinear eigenvalue systems and dense matrix operations. Furthermore, optimization techniques based on precomputations within the correction space render the total computational workload nearly independent of the number of self-consistent field iterations.
This approach significantly accelerates the solution process of the Hartree--Fock equation, effectively mitigating the traditional exponential scaling demands on computational resources while maintaining precision.


\end{abstract}

\section{Introduction}
The Hartree--Fock equation \cite{riginal1,original2} plays an important role in quantum physics, condensed matter physics, and quantum chemistry. It is used to address multielectron systems, particularly in determining the electronic structure of atoms, molecules, and condensed matter. The physical significance of the equation lies in its approach to approximately calculating the ground-state energy and wavefunction of a multielectron system by iteratively solving for the electron wavefunction and electron density. At its core idea, the method represents the multielectron wavefunction as a Slater determinant of single-electron wavefunctions.
With the rapid development of computer technology, the Hartree--Fock method has become widely used and extensively studied. Although the method has limitations, such as its inability to handle strongly correlated electron systems, it remains a fundamental tool for studying the electronic structures of multielectron systems. In particular, for simple atoms, molecules, and condensed matter, the Hartree--Fock method provides accurate results. It also serves as the foundation for most \emph{ab initio} wavefunction-based methods used to solve the electronic Schr\"{o}dinger equation.
The method also shares considerable similarities with density functional theory (DFT) approaches \cite{dft}, as the Kohn--Sham equations \cite{ks} differ only by the inclusion of an additional exchange-correlation potential. Under the local density approximation (LDA) \cite{ks} and the generalized gradient approximation (GGA) \cite{gga}, Kohn--Sham equations avoid the computation of the exact exchange potential required in Hartree--Fock equations. However, many modern functionals are hybrids that include a portion of the exact exchange functional. Solving Kohn--Sham equations with hybrid functionals \cite{hybrid} requires a computational effort similar to that of solving the Hartree--Fock equation, in addition to the numerical integration of the exchange-correlation potential. Therefore, designing an efficient numerical algorithm for the Hartree--Fock equation is also foundational for solving hybrid DFT.

There are many challenges in solving the Hartree--Fock equation numerically. First, it is a complex nonlinear system that describes interactions among electrons in a multielectron system, including both exchange and Coulomb repulsion interactions. This complexity necessitates the use of iterative numerical methods to obtain approximate solutions. Second, as the number of electrons in the system increases, the dimensionality of the Hartree--Fock equation grows rapidly, resulting in a sharp increase in computational and storage demands. Solving the Hartree--Fock equation for large-scale systems becomes extremely time-consuming, often exceeding available computational resources. Therefore, reducing computational cost and storage requirements while maintaining solution accuracy remains a key challenge. Additionally, the numerical solution process is influenced by the stability and convergence of the algorithm. Hence, developing efficient numerical methods capable of accurately simulating large-scale systems continues to be a major challenge in computational quantum physics and chemistry.

After discretization, the exchange interaction in the Hartree--Fock equation leads to matrices with many nonzero elements (i.e., dense matrices). These require more data to be processed during computation, thereby reducing
computational efficiency. In large-scale simulations, the multiplication and storage of dense matrices significantly increase computational costs. Consequently, most current methods for solving the Hartree--Fock equation rely on small basis set approaches, such as the local basis set method \cite{guass1,guass2,guass3,guass4,guass5}. Meanwhile, the plane-wave method is also widely used \cite{plane1,plane2,plane3,plane4,plane5,plane6,plane7,plane8}.
The finite element method (FEM) generally requires more degrees of freedom than the local basis set and plane-wave methods, rendering it extremely difficult to apply to the Hartree--Fock equation. However, FEM is particularly valuable when high-accuracy calculations are needed. By introducing more basis functions, FEM improves computational precision in describing the electronic structure of a system. Pask et al. \cite{pask} summarized the principles and advantages of FEM, noting its ability to provide solutions at predefined resolutions in real space using either periodic or nonperiodic wavefunctions. Nonperiodic functions are particularly useful for calculating the electronic properties of structures such as bent nanotubes or isolated molecules.
However, FEM results in a significant increase in computational workload and resource demand, particularly for large-scale systems. Thus, applying FEM to the Hartree--Fock equation remains relatively difficult. Developing efficient numerical algorithms tailored for FEM is therefore both important and challenging.

Previous efforts have explored integrating FEM with the Hartree--Fock method. Flores et al. \cite{flores1,flores2} used higher-order polynomials as basis functions to solve two-dimensional atomic Hartree--Fock equations. Heinemann et al. \cite{radial1} achieved high accuracy for light atoms and diatomic molecules using up to sixth-order Lagrange polynomials and prolate spheroidal coordinates, which significantly reduced computational cost but limited the range of applicable problems. In another study \cite{he14}, spin polarization was examined using the Hartree--Fock--Slater method. Lavor et al. \cite{lavor} proposed a global optimization approach using nonconvex polynomials as a variant of FEM. Sundholm et al. \cite{sun1,sun2} investigated finite-element multiconfiguration Hartree--Fock calculations for various atoms and compounds. Ackermann and Roitzsch \cite{ack1,ack2} applied a multigrid adaptive FEM in cylindrical coordinates for symmetric single-electron systems. Quiney et al. \cite{qui} compared finite basis set and FEM approaches for solving the relativistic Dirac equation. White et al. \cite{white} used second-order orthonormal shape functions on a uniform mesh for one- and two-electron systems. Although two-dimensional FEM-based Hartree--Fock calculations have been widely studied, the three-dimensional case is more complex and less explored.
In \cite{fem1}, a finite element Hartree--Fock method was introduced for polyatomic molecules using an accelerated exact exchange calculation with auxiliary local fitting. A divide-and-conquer method was also employed for parallelization and reduced computational scaling in complex systems. In \cite{nwch4}, Braun introduced a three-dimensional FEM approach for solving closed-shell Hartree--Fock equations in small molecules, using a factorized Green's function to solve both three-dimensional Schr\"{o}dinger-type and Poisson equations. Additional results are available in \cite{other4,other2,wp1,wp2,other5,other1,other3} and the references cited therein. As the Hartree--Fock equation is studied across disciplines such as computational mathematics, physics, and chemistry, only a portion of the relevant literature can be cited here.

To improve the accuracy and efficiency of numerical simulations of the three-dimensional Hartree--Fock equation, we adopt the multilevel correction method described in \cite{xiaohong,HuXieXu2,siam}, combined with an efficient implementation strategy. We also incorporate an adaptive refinement technique to address singularities \cite{babuska,cason}. Using the multilevel correction method, we solve a set of linearized boundary value problems and refine their approximate solutions by solving small-scale Hartree--Fock equations in a low-dimensional correction space. This correction space comprises a coarse space and the solution to the linearized boundary value problem. This approach avoids solving the large-scale equation directly and significantly improves computational efficiency.
In the correction space, we perform self-consistent field (SCF) iterations to solve the small-scale Hartree--Fock equation. To further enhance efficiency, we introduce two acceleration strategies. First, we construct separate correction spaces for each wavefunction independently. This ensures that the dimensionality of each correction space remains low and independent of the number of wavefunctions being solved. Moreover, the independent structure of these spaces facilitates parallel computation.
Second, because the coarse space that forms part of the correction space remains fixed, the approximate wavefunctions that serve as basis functions for the correction space are retained after each adaptive refinement. A large portion of the computational workload can be preprocessed using tensor computation techniques. Subsequent iterations then require only simple matrix and tensor operations to complete the solution process. This strategy renders the computational cost of the algorithm independent of the number of SCF iterations.
By employing multilevel correction techniques, we avoid large-scale nonlinear eigenvalue problems and dense matrix operations. Therefore, the multilevel correction algorithm significantly accelerates the solution of the Hartree--Fock equation.

The remainder of this paper is organized as follows: In Section 2, we review fundamental concepts and the finite element method for the
Hartree--Fock equation. Section 3 discusses the construction of the corresponding multilevel correction type of adaptive finite element method. In Section 4, we present an efficient implementation of the correction step within the correction space. Section 5 introduces numerical experiments to demonstrate the effectiveness of the proposed method. Finally, Section 6 provides concluding remarks.


\section{Finite element method for Hartree--Fock equation}
In this study, we consider a molecular system comprising $M$ nuclei with charges ${Z_1,\cdots,Z_M}$ and locations ${R_1, \cdots, R_M}$ with
$N$ electrons.
The energy of the considered molecular system is described by:
\begin{eqnarray}
 E_{HF} &=& \int_{\mathbb R^3} \left( \frac{1}{2}\sum_{\ell=1}^N|\nabla \phi_\ell |^2-\sum_{k=1}^{M}\frac{Z_k}{|x-R_k|}\rho(x) \right)dx
 +\frac{1}{2}\int_{\mathbb R^3}\int_{\mathbb R^3}\frac{\rho(x)\rho(y)}{|x-y|}dxdy \nonumber\\
 &&
 -\frac{1}{2}\int_{\mathbb R^3}\int_{\mathbb R^3}\frac{|P(x,y)|^2}{|x-y|}dxdy,
\end{eqnarray}
where the density matrix $P(x,y)=\sum_{\ell=1}^N \phi_\ell(x)\phi_\ell^*(y)$,
the density function $\rho(x)=P(x,x)=\sum_{\ell=1}^N|\phi_\ell(x)|^2$.

The corresponding Euler--Lagrange equation for the minimum Hartree--Fock energy is just the Hartree--Fock equation:
\begin{eqnarray}\label{defHartree-Focke}
 -\frac{1}{2}\Delta \phi_\ell +V_{\rm ext} \phi_\ell
 +V_{\rm Har}(\rho) \phi_\ell +V_{\rm x}(P)\phi_\ell=\lambda_\ell \phi_\ell,\quad \ell=1,\cdots,N,
\end{eqnarray}
where the external potential $V_{\rm ext}$, the Hartree potential $V_{\rm Har}(\rho)$ and the exchange potential $V_{\rm x}(P)$ are defined by:
\begin{eqnarray}
&&V_{\rm ext} = -\sum_{k=1}^{M}\frac{Z_k}{|x-R_k|},\label{defext}\\
&&V_{\rm Har}(\rho) = \int_{\mathbb R^3}\frac{\rho(y)}{|x-y|}dy,\label{defhar}\\
&&V_{\rm x}(P)\phi_\ell = -\int_{\mathbb R^3}\frac{P(x,y)}{|x-y|}\phi_\ell(y)dy.\label{defx}
\end{eqnarray}

As wavefunctions decay exponentially to zero, we assume that the numerical simulations are implemented in the bounded domain $\Omega$.
To apply the FEM, we provide a variational form of the Hartree--Fock equation (\ref{defHartree-Focke}) as follows:
\begin{eqnarray}
 \frac{1}{2}(\nabla \phi_\ell,\nabla \psi) +(V_{\rm ext} \phi_\ell, \psi)
 +(V_{\rm Har}(\rho) \phi_\ell, \psi) +(V_{\rm x}(P)\phi_\ell, \psi)=(\lambda_\ell \phi_\ell, \psi),\quad \forall \psi\in H_0^1(\Omega),
\end{eqnarray}
where $(\phi,\psi)=\int_\Omega \phi\psi^*d\Omega$ denotes the standard $L^2$-inner product.

Next, we generate a decomposition $\mathcal T_h$ in computing domain $\Omega$. Using mesh $\mathcal T_h$,
we construct a finite element space $V_h$ which is composed of piecewise polynomials.
Then, the FEM is used to solve the following discrete Hartree--Fock equation:
\begin{eqnarray*}
 \frac{1}{2}(\nabla \phi_{\ell,h},\nabla \psi_h) +(V_{\rm ext} \phi_{\ell,h}, \psi_h)
 +(V_{\rm Har}(\rho_h) \phi_{\ell,h}, \psi_h) +(V_{\rm x}(P_h)\phi_{\ell,h}, \psi_h)=(\lambda_{\ell,h} \phi_{\ell,h}, \psi_h),\quad \forall \psi_h\in V_h.
\end{eqnarray*}

Because the discrete Hartree--Fock equation is a nonlinear eigenvalue problem, the SCF iteration should be adopted
to derive the approximate solutions gradually. Hence, a linearized eigenvalue problem should be solved at each iteration step. The standard SCF iteration
is presented in Algorithm \ref{Multilevel Correction1} as follows:
\begin{algorithm}\label{Multilevel Correction1} SCF Iteration for Hartree--Fock equation
\begin{enumerate}
\item Given an initial value $\phi_{\ell,h}^{(1)}, \ \ell=1,2,\cdots,N$, and a tolerance $TOL$.
\item Set $p=1$.
\item Derive the new approximate wavefunctions $\phi_{\ell,h}^{(p+1)}, \ \ell=1,2,\cdots,N$ by solving the following linearized eigenvalue problems:
\begin{eqnarray}
 && \frac{1}{2}(\nabla \phi_{\ell,h}^{(p+1)},\nabla \psi_h) +(V_{\rm ext} \phi_{\ell,h}^{(p+1)}, \psi_h)
 +(V_{\rm Har}(\rho_h^{(p)}) \phi_{\ell,h}^{(p+1)}, \psi_h) +(V_{\rm x}(P_h^{(p)})\phi_{\ell,h}^{(p+1)}, \psi_h) \nonumber\\
 &=& (\lambda_{\ell,h}^{(p+1)} \phi_{\ell,h}^{(p+1)}, \psi_h),\quad \forall \psi_h\in V_h.
\end{eqnarray}

\item if $|E^{p+1}-E^{p}|/|E^{p+1}|\leq TOL$, stop. Else set $p=p+1$ and go to Step 3.
\end{enumerate}
\end{algorithm}

For the Hartree potential $V_{\rm Har}(\rho_h)$, we can derive an approximation by solving the following Poisson equation:
\begin{eqnarray}\label{harbd}
-\Delta V_{\rm Har}(\rho_h) = 4\pi\rho_{h}(x).
\end{eqnarray}
Since the wavefunction of the Hartree--Fock equation decreases exponentially,
zero Dirichlet boundary conditions can be adopted when the computational domain $\Omega$ is sufficiently large.
However, since the solution of equation (\ref{harbd}) decays linearly, simply using the zero Dirichlet boundary condition introduces a large truncation error at the boundary.
In order to guarantee the accuracy, we approximate boundary values using multipole expansion \cite{Baogang}.

For the exchange potential, based on the definition of density matrix, we obtain
\begin{eqnarray}\label{Avddef}
\big(V_{\rm x}(P)\psi_j, \psi_i\big) &=& -\big(\int_{\Omega}\frac{\rho(x,y)}{|x-y|}\psi_j(y)dy, \psi_i(x)\big) \nonumber\\
&=&  -\big( \int_{\Omega}\frac{\sum_{\ell=1}^N \phi_\ell(x)\phi_\ell^*(y)}{|x-y|}\psi_j(y)dy , \psi_i(x)\big) \nonumber\\
&=& -\sum_{\ell=1}^N \int_{\Omega} \bigg( \int_{\Omega}\frac{ \phi_\ell^*(y)\psi_j(y)}{|x-y|} dy \bigg) \phi_\ell(x)\psi_i(x)dx.
\end{eqnarray}
Thus, a dense matrix will be derived by discretizing the exchange potential. Meanwhile, to discretize (\ref{Avddef}), we need
to first solve the functions $\int_{\Omega}\frac{ \phi_\ell^*(y)\psi_j(y)}{|x-y|} dy$, which requires solving $NN_h$ Poisson equations with $N_h={\rm dim}(V_h)$.
Thus, it is quite a time-consuming step to discretize the exchange potential.

\section{Multilevel correction type of adaptive finite element method for Hartree--Fock equation}
The numerical solution of the Hartree--Fock equation investigated in this study is extremely challenging.
First, it is a nonlinear eigenvalue model that must be solved using SCF iterations and mixing schemes \cite{mix1,mix2}.
Meanwhile, each iteration step involves solving a linear eigenvalue problem.
As we know, it is quite difficult to solve large-scale eigenvalue problems, and the computational workload increases exponentially.
Next, owing to the presence of the exchange potential, the matrix obtained after discretizing the eigenvalue model becomes dense.
This significantly increases the difficulty of addressing storage-related challenges.
When the numerical algorithm involves a large amount of basis functions, such as with the FEM,
numerical simulations become highly demanding. Thus, to date,
there are very few finite element algorithms available for solving the Hartree--Fock equation in three-dimensional space.

In order to improve the efficiency and accuracy of the numerical simulation, we design a multilevel correction type of adaptive algorithm for solving the Hartree--Fock equation in this paper.
We first introduce the multilevel correction technique to decompose the Hartree--Fock equation to reduce the computational difficulty.
Subsequently, the adaptive refinement technique is combined to address the singularities, thereby forming a complete multilevel correction adaptive algorithm.
The idea of our algorithm is to solve the linearized boundary value problem in the adaptive finite element space to get the approximate solution of the wavefunction.
The exchange potential is placed on the righthand side of the linearized equation to reduce the computational complexity.
Thus, this linearized boundary value problem requires significantly less computational effort compared with a full nonlinear eigenvalue problem.
However, using a boundary value problem to approximate a nonlinear eigenvalue problem inevitably introduces errors.
In order to further improve the accuracy, we construct a low-dimensional correction space, which consists of a low-dimensional finite element space
and the approximate solution of the linearized boundary value problem. So this correction space is highly accurate while maintaining a low dimensionality.
Then we solve a small-scale Hartree--Fock equation in this correction space to get the final approximate solution of the current adaptive finite element space.
Therefore, the difficult parts of solving the Hartree--Fock equation are all confined to a low-dimensional and high-precision space.

For each wavefunction, we construct an individual correction space based on the corresponding solution of the linearized boundary value problem, and correct each wavefunction separately. This approach avoids high-dimensional inner product calculations that arise in global solution procedures and preserves the orthogonality of different wavefunctions. Additionally, treating each wavefunction with a separate correction process is well-suited for parallel computation.

The algorithm operates continuously on a multilevel mesh sequence until the prescribed tolerance is achieved.
Assume that we have obtained an approximate solution $(\lambda_{\ell,h_{k}}, \phi_{\ell,h_{k}})$
in the space $V_{h_{k}}$. The detailed process in space $V_{h_{k+1}}$ to derive a new approximation $(\lambda_{\ell,h_{k+1}}, \phi_{\ell,h_{k+1}})$ is described as follows:
\begin{algorithm}\label{Multilevel Correction2} One correction step for Hartree--Fock equation
\begin{enumerate}
\item Solve the following linearized boundary value problem: 
\begin{eqnarray}\label{Aux_Linear_Problem}
 && \frac{1}{2}(\nabla \widetilde{\phi}_{\ell,h_{k+1}}, \nabla \psi_{h_{k+1}}) +(V_{\rm ext} \widetilde{\phi}_{\ell,h_{k+1}}, \psi_{h_{k+1}})
 +(V_{\rm Har}(\rho_{h_{k}}) \widetilde{\phi}_{\ell,h_{k+1}}, \psi_{h_{k+1}})  \nonumber\\
 &=& (\lambda_{\ell,h_{k}}\phi_{\ell,h_{k}}, \psi_{h_{k+1}})-(V_{\rm x}(P_{h_{k}})\phi_{\ell,h_{k}}, \psi_{h_{k+1}}),\quad \forall \psi_{h_{k+1}}\in V_{h_{k+1}}.
\end{eqnarray}
\item Define a correction space $V_{H,h_{k+1}}=V_H+{\rm span}\{\widetilde{\phi}_{\ell,h_{k+1}}\}$ and
solve the following small-scale Hartree--Fock equation:
\begin{eqnarray}\label{Nonlinear_Eig_Hh}
 && \frac{1}{2}(\nabla {\phi}_{\ell,h_{k+1}}, \nabla \psi_{H,h_{k+1}}) +(V_{\rm ext} {\phi}_{\ell,h_{k+1}}, \psi_{H,h_{k+1}})\nonumber\\
 && +(V_{\rm Har}(\rho_{h_{k+1}}) {\phi}_{\ell,h_{k+1}}, \psi_{H,h_{k+1}}) +(V_{\rm x}(P_{h_{k+1}}){\phi}_{\ell,h_{k+1}}, \psi_{H,h_{k+1}}) \nonumber\\
 &=& (\lambda_{\ell,h_{k+1}}\phi_{\ell,h_{k+1}}, \psi_{H,h_{k+1}}),\quad \forall \psi_{h_{k+1}}\in V_{H,h_{k+1}}.
\end{eqnarray}
\end{enumerate}
\end{algorithm}
For simplicity, we summarize the two steps above using the following expressions:
\begin{eqnarray}\label{Nonlinear_Eig_Hhsymbol}
(\lambda_{\ell,h_{k+1}},\phi_{\ell,h_{k+1}})=HFCorrection(\lambda_{\ell,h_{k}},\phi_{\ell,h_{k}},V_H,V_{h_{k+1}}).
\end{eqnarray}

As mentioned, the novel algorithm is executed on a multilevel mesh sequence. Therefore, the mesh must be gradually refined to generate this sequence. Owing to the singularity of the Hartree--Fock equation, the adaptive finite element method (AFEM) is a competitive choice.
Based on an \emph{a posteriori} error estimator, AFEM selects mesh elements with large errors for further refinement. After generating the new mesh, we solve the Hartree--Fock equation and obtain a new approximate solution. We then begin another loop by computing the \emph{a posteriori} error estimator and selecting mesh elements for additional refinement. This iterative process can be described as follows:
\begin{center}
$\cdots$\bf{Solve} $\rightarrow$ \bf{Estimate}
$\rightarrow$ \bf{Mark} $\rightarrow$ \bf{Refine}$\cdots$.
\end{center}
To obtain $\mathcal T_{h_{k+1}}$ from $\mathcal T_{h_k}$, we first solve the discrete equation on $\mathcal T_{h_k}$ to derive an approximate solution, and then compute the \emph{a posteriori} error estimator for each mesh element. Next, we mark the elements with the largest errors and refine them such that the new mesh $\mathcal T_{h_{k+1}}$ remains shape-regular and conforming.

In our simulation, the residual type \emph{a posteriori}
  error estimation is employed to generate the error indicator. First,
we construct the element residual ${\mathcal R}_{T}(\{\lambda_{\ell,h_k},\phi_{\ell,h_k}\}_{\ell=1}^N)$ and the jump
residual ${\mathcal J}_e(\{\phi_{\ell,h_k}\}_{\ell=1}^N)$ for the eigenpair approximation
$(\lambda_{\ell,h_k},\phi_{\ell,h_k}), \ \ell=1,\cdots,N$ as follows:
\begin{eqnarray}
&&{\mathcal R}_{T}(\{\lambda_{\ell,h_k},\phi_{\ell,h_k}\}_{\ell=1}^N) := \big(-\frac{1}{2}\Delta \phi_{\ell,h_k} +V_{\rm ext} \phi_{\ell,h_k}
 +V_{\rm Har}(\rho_{h_k}) \phi_{\ell,h_k} +V_{\rm x}(P_{h_k})\phi_{\ell,h_k}   \nonumber\\
&& \ \ \ \ \ \ \ \ \ \ \  \ \ \ \ \ \ \ \ \ \ \  \ \ \ \ \ \ \ \ \ \ \ \ \ -\lambda_{\ell,h_k}\phi_{\ell,h_k}\big)_{\ell=1}^{N},
\ \ \ \ \ \  \ \ \ \ \ \ \ \ \ \ \ \ \ \ \ \ \ \ \ \ \ \ \ \ \  \text{in } T\in \mathcal T_{h_k}, \\
&&{\mathcal J}_e(\{\phi_{\ell,h_k}\}_{\ell=1}^N) := \big( \frac{1}{2}\nabla \phi_{\ell,h_k}|_{T^+}\cdot \nu^+
+\frac{1}{2}\nabla \phi_{\ell,h_k}|_{T^-}\cdot \nu^- \big)_{\ell=1}^{N}, \ \  \ \ \ \ \ \ \ \ \text{on } e\in \mathcal E_{h_k},
\end{eqnarray}
where $e$ is the common side of elements $T^+$ and $T^-$ with the unit outward normals $\nu^+$ and
$\nu^-$, $\mathcal E_{h_k}$ denotes the set of interior faces of $\mathcal T_{h_k}$.
For mesh element $T\in\mathcal T_{h_k}$, we define the local error indicator ${\eta}_{h_k}^2(\{\lambda_{\ell,h_k},\phi_{\ell,h_k}\}_{\ell=1}^N,T)$ by
\begin{eqnarray*}
{\eta}_{h_k}^2(\{\lambda_{\ell,h_k},\phi_{\ell,h_k}\}_{\ell=1}^N,T):=h_T^2\| {\mathcal R}_{T}(\{\lambda_{\ell,h_k},\phi_{\ell,h_k}\}_{\ell=1}^N)\|_{0,T}^2+\sum_{e\in\mathcal E_{h_k},e
\subset \partial T}h_e\| {\mathcal J}_{e}(\{\phi_{\ell,h_k}\}_{\ell=1}^N)\|_{0,e}^2.
\end{eqnarray*}
Given a subset $\widehat\Omega \subset\Omega$, we define the error estimate ${\eta}_{h_k}^2(\{\lambda_{\ell,h_k},\phi_{\ell,h_k}\}_{\ell=1}^N,\widehat\Omega)$ by
\begin{eqnarray}\label{def of eta and osc}
{\eta}_{h_k}^2(\{\lambda_{\ell,h_k},\phi_{\ell,h_k}\}_{\ell=1}^N,\widehat\Omega)=\sum_{T\in \mathcal T_{h_k},T\subset \widehat\Omega}{\eta}_{h_k}^2(\{\lambda_{\ell,h_k},\phi_{\ell,h_k}\}_{\ell=1}^N,T).
\end{eqnarray}
Based on the error indicator (\ref{def of eta and osc}), we use the
D\"{o}rfler's marking strategy \cite{dorfler} described in Algorithm \ref{dorfler} to mark certain elements for local refinement.

\begin{algorithm}\label{dorfler}D\"{o}rfler's Marking Strategy
\begin{enumerate}
\item
Given a parameter $\theta\in(0,1)$.
\item
Construct a minimal subset ${\mathcal M}_{h_k}$ from  $\mathcal T_{h_k}$ by selecting some elements in $\mathcal T_{h_k}$ such that
\begin{eqnarray*}
\sum_{T\in \mathcal M_{h_k}}{\eta}_{h_k}^2(\{\lambda_{\ell,h_k},\phi_{\ell,h_k}\}_{\ell=1}^N,T)\geq \theta {\eta}_{h_k}^2(\{\lambda_{\ell,h_k},\phi_{\ell,h_k}\}_{\ell=1}^N,\Omega).
\end{eqnarray*}
\item
 Mark all the elements in ${\mathcal M}_{h_k}.$
\end{enumerate}
\end{algorithm}
Based on the adaptive refinement scheme described above, and the one level
correction step defined by Algorithm \ref{Multilevel Correction2}, the multilevel
correction type of AFEM for Hartree--Fock equation is given below:

\begin{algorithm}\label{multilevelAFEM}Multilevel correction type of AFEM for Hartree--Fock equation
\begin{enumerate}
\item Generate a coarse mesh $\mathcal T_H$ on $\Omega$ and construct the corresponding coarse space $V_H$.
Start with an initial mesh $\mathcal T_{h_1}$  obtained by refining $\mathcal T_H$  several times.
Subsequently, establish the initial space $V_{h_1}$ on $\mathcal T_{h_1}$  and solve the Hartree--Fock equation:
\begin{eqnarray}\label{f1Hartree--Fock}
 && \frac{1}{2}(\nabla {\phi}_{\ell,h_1}, \nabla \psi_{h_1}) +(V_{\rm ext} {\phi}_{\ell,h_1}, \psi_{h_1}) 
  +(V_{\rm Har}(\rho_{h_{1}}) {\phi}_{\ell,h_1}, \psi_{h_1}) +(V_{\rm x}(P_{h_{1}}){\phi}_{\ell,h_1}, \psi_{h_1}) \nonumber\\
 &=& (\lambda_{\ell,h_{1}}\phi_{\ell,h_{1}}, \psi_{h_1}),\ \ \ell=1,\cdots,N,\quad \forall \psi_{h_1}\in V_{h_1}.
\end{eqnarray}
\item
Set $k=1$.
\item
Compute the local error indicator ${\eta}_{h_k}(\{\lambda_{\ell,h_k},\phi_{\ell,h_k}\}_{\ell=1}^N,T)$ for each mesh element $T\in\mathcal T_{h_k}$.
\item
Construct $\mathcal M_{h_k}\subset \mathcal T_{h_k}$ by D\"{o}rfler's marking strategy and
then refine $\mathcal T_{h_k}$ to get a new mesh $\mathcal T_{h_{k+1}}$.
\item For $\ell=1,\cdots,N$, perform the following correction step to obtain $(\lambda_{\ell,h_{k+1}},\phi_{\ell,h_{k+1}})$:
\begin{eqnarray}
(\lambda_{\ell,h_{k+1}},\phi_{\ell,h_{k+1}})=HFCorrection(\lambda_{\ell,h_{k}},\phi_{\ell,h_{k}},V_H,V_{h_{k+1}}).
\end{eqnarray}
\item
Let $k=k+1$ and go to Step 3.
\end{enumerate}
\end{algorithm}

\section{Efficient implementation of correction step}
This section introduces the efficient strategy to implement the multilevel correction type of AFEM defined in the previous section.
In each adaptive finite element space, we only need to solve some linearized boundary value problems without nonlinear iteration.
Additionally, to avoid generating dense matrices,
the exchange potential is placed on the righthand side of the equation. Thus, the computation time and storage space
are significantly smaller compared to those required for solving the large-scale Hartree--Fock equation.

For the correction step defined in $V_{H,h_k}$, the SCF iteration is adopted to solve the small-scale Hartree--Fock equation.
Then in each iteration step, we need to solve a small-scale linear dense eigenvalue problem. Because the dimension of the correction space ($N_H+1$)
is small compared to that of the adaptive finite element space ($N_{h_k}$), thus the computational efficiency and storage requirement
are significantly improved relative to those of the large-scale dense eigenvalue problem.
In the correction space, the solution of the linearized boundary value problem is added as a basis function,
thus the integral calculations associated with this basis function should be performed on the adaptive mesh to maintain the accuracy.
Additionally, complex molecular systems often require numerous SCF iterations, with each iteration repeating the solution process.
Hence, for the multilevel correction type of AFEM,
the most time-consuming step is the correction process defined in the low-dimensional correction space if we adopt the traditional SCF scheme described in Algorithm \ref{Multilevel Correction1}.

To improve the efficiency of solving the small-scale Hartree--Fock equation,  we next present an efficient strategy for implementing the correction step.
As the coarse space $V_H$ remains unchanged throughout the algorithm and the newly added basis function remains fixed during the SCF iteration,
we can perform the invariant computations outside the SCF iteration.
Hence, each iteration step requires only a few simple matrix operations.
In addition, most steps in the computation process support parallelization.
For example, the linearized boundary problems and correction spaces are all constructed separately for each wavefunction.

Next, we will use the $\ell$-th wavefunction as an example to illustrate the implementation process of the correction step, which is the same for all the other wavefunctions.
For the description of implementing technique here, let  $\{\psi_{j,H}\}_{1\leq j\leq N_H}$ denote the Lagrange basis function for the coarse finite element space $V_H$,
and $\widetilde \phi_{\ell,h}$ denote the added basis function in the correction space $V_{H,h}$.
Then the correction space for the $\ell$-th wavefunction can be described by $V_{H,h}=\{\widetilde \phi_{\ell,h},\psi_{1,H},\cdots,\psi_{N_H,H}\}$.

%
In order to show the main idea here, we consider the matrix version of eigenvalue problem (\ref{Nonlinear_Eig_Hh}).
Since a new basis function is added, the corresponding matrix should have one more row and one more column.
The detailed form can be described as follows:
\begin{equation}\label{Eigenvalue_H_h}
\left(
\begin{array}{cc}
A_H & b_{Hh}\\
b_{Hh}^H&\beta
\end{array}
\right)
\left(
\begin{array}{c}
C_{H} \\
\theta
\end{array}
\right)
=\lambda_{\ell,h}\left(
\begin{array}{cc}
M_H & c_{Hh}\\
c_{Hh}^H&\gamma
\end{array}
\right)
\left(
\begin{array}{c}
C_H\\
\theta
\end{array}
\right),
\end{equation}
where dim$(A_H)$=dim$(M_H)$=$N_H\times N_H$, dim$(b_{Hh})$=dim$(c_{Hh})$=$N_H\times 1$, dim$(C_H)=N_H\times 1$ and dim$(\beta)$= dim$(\gamma)$= dim$(\theta)=1$.


Based on the Hartree--Fock equation (\ref{Nonlinear_Eig_Hh}), the matrix
$M_H$ has the following form:
\begin{eqnarray}\label{massmatrix}
(M_H)_{i,j}=\int_\Omega\psi_{j,H} \psi_{i,H}dx,\ \ (c_{Hh})_{i}=\int_\Omega\widetilde \phi_{\ell,h} \psi_{i,H}dx,\ \ \gamma=\int_\Omega|\widetilde \phi_{\ell,h}|^2dx.
\end{eqnarray}
It is obvious that the matrix $M_H$ will not change throughout the algorithm.
The vector $c_{Hh}$ and the scalar $\gamma$ rely on the basis function $\widetilde \phi_{\ell,h}$.
As the basis function $\widetilde \phi_{\ell,h}$ remains unchanged in the correction space,
thus the vector $c_{Hh}$ and the scalar $\gamma$ will not change during the SCF iteration in the correction space.

Because the wavefucntion will be renewed at each iteration step,
the stiffness matrix $A_H$, the vector $b_{Hh}$ and the scalar $\beta$ will change during the SCF iteration process.
Then it is required to consider the efficient implementation to update
the matrix $A_H$, the vector $b_{Hh}$ and the scalar $\beta$.

%
%
Assume we have a given initial wavefunction $\phi_{\ell,Hh}\in V_{H,h}$ for the SCF iteration involved in the correction space.
The approximate wavefunction $\phi_{\ell,Hh}$ can be expressed as $\phi_{\ell,Hh}=\widetilde \theta \widetilde\phi_{\ell,h}+ \sum_{m=1}^{N_H}c_{m,H}\psi_{m,H}$.
As $\phi_{\ell,Hh}$ is the given initial value, the coefficients $c_{1,H},\cdots,c_{m,H},\widetilde\theta$ are known,
which is just the solution of (\ref{Eigenvalue_H_h}) in the
previous iteration step.
In the subsequent iterative solving process, we will continually update the value of $\phi_{\ell,Hh}$ (i.e., update the coefficients $c_{1,H},\cdots,c_{m,H},\widetilde\theta$).
When we correct the $\ell$-the wavefunction $\phi_{\ell,Hh}$ in its own correction space, the other wavefunctions are fixed as $\widetilde\phi_{s,h}$ (i.e. $\phi_{s,Hh}=\widetilde\phi_{s,h}$ if $s\neq \ell$)
during the computation of the density function and density matrix.

From the definitions of the correction space $V_{H,h}$ and the eigenvalue problem (\ref{Eigenvalue_H_h}), the matrix
$A_H$ has the following expansion:
\begin{eqnarray}\label{Adefn1}
&&(A_H)_{i,j}\nonumber\\
&=&\frac{1}{2}\int_\Omega\nabla\psi_{j,H}(x)\nabla \psi_{i,H}(x)dx +\int_\Omega V_{\rm ext}\psi_{j,H}(x)\psi_{i,H}(x)dx \nonumber\\
&&+ \int_\Omega V_{\rm Har}(\rho_{Hh})\psi_{j,H}(x)\psi_{i,H}(x)dx +\int_{\Omega}(V_{\rm x}(P_{Hh})\psi_{j,H})\psi_{i,H}(x)dx \nonumber\\
&=&\frac{1}{2}\int_\Omega\nabla\psi_{j,H}(x)\nabla \psi_{i,H}(x)dx +\int_\Omega V_{\rm ext}\psi_{j,H}(x)\psi_{i,H}(x)dx \nonumber\\
&&+ \sum_{s=1}^N\int_\Omega\int_\Omega \frac{ \phi_{s,Hh}(y)\phi_{s,Hh}^*(y)}{|x-y|}dy\psi_{j,H}(x)\psi_{i,H}(x)dx \nonumber\\
&&+ \sum_{s=1}^N\int_\Omega\int_\Omega \frac{ \phi_{s,Hh}^*(y)\psi_{j,H}(y)}{|x-y|}dy\phi_{s,Hh}(x)\psi_{i,H}(x)dx.\ \ \ \ \ \ \
\end{eqnarray}

Let us denote
\begin{eqnarray}
(A_{\rm Kin})_{i,j}&=&\frac{1}{2}\int_\Omega\nabla\psi_{j,H}(x)\nabla \psi_{i,H}(x)dx +\int_\Omega V_{\rm ext}\psi_{j,H}(x)\psi_{i,H}(x)dx, \label{akin}\\
(A_{\rm Har,s})_{i,j}&=&\int_\Omega\int_\Omega \frac{\phi_{s,Hh}(y)\phi_{s,Hh}^*(y)}{|x-y|}dy\psi_{j,H}(x)\psi_{i,H}(x)dx, \label{ahar}\\
(A_{\rm x,s})_{i,j}&=&\int_\Omega\int_\Omega \frac{ \phi_{s,Hh}^*(y)\psi_{j,H}(y)}{|x-y|}dy\phi_{s,Hh}(x)\psi_{i,H}(x)dx.\label{ax}
\end{eqnarray}
The matrix $A_{\rm Kin}$ remains unchanged throughout the algorithm. Thus, we only need to calculate $A_{\rm Kin}$ one time at the beginning of the algorithm.
When $s\neq\ell$,   we have
\begin{eqnarray}\label{Adefsl}
&&(A_{\rm Har,s})_{i,j} + (A_{\rm x,s})_{i,j} \nonumber\\
&=&\int_\Omega\int_\Omega \frac{| \widetilde\phi_{s,h}(y)|^2}{|x-y|}dy\psi_{j,H}(x)\psi_{i,H}(x)dx \nonumber\\
&& + \int_\Omega \int_\Omega \frac{\widetilde\phi_{s,h}^*(y)\psi_{j,H}(y)}{|x-y|}dy  \widetilde\phi_{s,h}(x)\psi_{i,H}(x) dx.  \ \ \ \ \ \
\end{eqnarray}
When $s=\ell$,   we have
\begin{eqnarray}\label{Adef}
&&(A_{\rm Har,\ell})_{i,j} + (A_{\rm x,\ell})_{i,j} \nonumber\\
&=&\int_\Omega\int_\Omega \frac{\Big|\widetilde \theta \widetilde\phi_{\ell,h}(y)+ \sum_{m=1}^{N_H}c_{m,H}\psi_{m,H}(y)\Big|^2}{|x-y|}dy\psi_{j,H}(x)\psi_{i,H}(x)dx \nonumber\\
&& + \int_\Omega \Big\{\int_\Omega \frac{\Big(\widetilde \theta \widetilde\phi_{\ell,h}(y)+ \sum_{m=1}^{N_H}c_{m,H}\psi_{m,H}(y)\Big)^*\psi_{j,H}(y)}{|x-y|}dy \nonumber\\
&&\ \ \ \ \ \ \ \  \Big(\widetilde \theta\widetilde\phi_{\ell,h}(x)+ \sum_{m=1}^{N_H}c_{m,H}\psi_{m,H}(x)\Big)\psi_{i,H}(x) \Big\}dx.  \ \ \ \ \ \
\end{eqnarray}
For simplicity, let us denote
\begin{eqnarray*}
\widetilde{\widetilde V}_{s,\ell}(x) &=& \int_\Omega\frac{\widetilde\phi_{s,h}(y)\widetilde\phi_{\ell,h}^*(y)}{|x-y|}dy,\nonumber\\
\widetilde V_{s,m}(x) &=& \int_\Omega\frac{\widetilde\phi_{s,h}(y)\psi_{m,H}(y)}{|x-y|}dy,\nonumber\\
V_{m,n}(x) &=& \int_\Omega\frac{\psi_{m,H}(y)\psi_{n,H}(y)}{|x-y|}dy.
\end{eqnarray*}
Then (\ref{Adefsl}) and (\ref{Adef}) can be rewritten as
\begin{eqnarray}\label{Adef2sl}
&&(A_{\rm Har,s})_{i,j} + (A_{\rm x,s})_{i,j} \nonumber\\
&=&\int_\Omega \widetilde{\widetilde V}_{s,s}(x)\psi_{j,H}(x)\psi_{i,H}(x)dx  + \int_\Omega \widetilde V_{s,j}^*(x)  \widetilde\phi_{s,h}(x)\psi_{i,H}(x) dx\nonumber\\
\end{eqnarray}
and
\begin{eqnarray}\label{Adef2}
&&(A_{\rm Har,\ell})_{i,j} + (A_{\rm x,\ell})_{i,j} \nonumber\\
&=&\int_\Omega \Big\{\Big(|\widetilde \theta|^2\widetilde{\widetilde V}_{\ell,\ell}(x) + \widetilde \theta\sum_{m=1}^{N_H} c_{m,H}^*\widetilde V_{\ell,m}(x)  + \widetilde \theta^*\sum_{m=1}^{N_H} c_{m,H}\widetilde V_{\ell,m}^*(x)
+ \sum_{m,n=1}^{N_H}c_{m,H}c_{n,H}^* V_{m,n}(x) \Big) \nonumber\\
&&\ \ \ \ \  \psi_{j,H}(x)\psi_{i,H}(x)\Big\}dx \nonumber\\
&& + \int_\Omega \Big(\widetilde \theta^* \widetilde V_{\ell,j}^*(x) +\sum_{m=1}^{N_H} c_{m,H}^*V_{j,m}(x)  \Big) \Big(\widetilde\theta\widetilde\phi_{\ell,h}(x)
+ \sum_{m=1}^{N_H}c_{m,H}\psi_{m,H}(x)\Big)\psi_{i,H}(x)dx.
\end{eqnarray}
To further simplify the expression, we define
\begin{eqnarray}
V_{m,n,j,i} &=& \int_\Omega V_{m,n}(x)\psi_{j,H}(x)\psi_{i,H}(x)dx, \label{vmnjk} \\
\widetilde V_{\ell,m,j,i} &=& \int_\Omega\widetilde V_{\ell,m}(x)\psi_{j,H}(x)\psi_{i,H}(x)dx,\label{lmkj}\\
\widetilde{\widetilde V}_{\ell,j,i} &=& \int_\Omega\widetilde{\widetilde V}_{\ell,\ell}(x)\psi_{j,H}(x)\psi_{i,H}(x)dx,\label{vlkj}\\
\widetilde W_{\ell,j,i} &=& \int_\Omega \widetilde V_{\ell,j}(x)\widetilde \phi_{\ell,h}(x)\psi_{i,H}(x)dx,\label{wlkj}\\
\widetilde W_{*,\ell,j,i} &=& \int_\Omega \widetilde V_{\ell,j}^*(x)\widetilde \phi_{\ell,h}(x)\psi_{i,H}(x)dx.\label{*wlkj}
\end{eqnarray}
Based on the definition of $\widetilde V_{\ell,m,j,i}$, we can also derive
$\widetilde V_{\ell,i,j,m} = \int_\Omega V_{j,m}(x)\widetilde\phi_{\ell,h}(x)\psi_{i,H}(x)dx.$

As we can see, $V_{m,n,j,i}$ remains unchanged during the whole algorithm. Thus, we can calculate $V_{m,n,j,i}$ at the beginning of the algorithm.
Because $\widetilde\phi_{\ell,h}$ remains unchange in each correction step, so $\widetilde V_{\ell,m,j,i}$, $\widetilde{\widetilde V}_{\ell,j,i}$,
$\widetilde W_{\ell,j,i}$ and $\widetilde W_{*,\ell,j,i}$ can be calculated before each correction step.

After calculating these quantities in advance, we can assemble $A_{\rm Har,\ell}$ and $A_{\rm x,\ell}$ by:
\begin{eqnarray}\label{Adef2sldd}
(A_{\rm Har,s})_{i,j} + (A_{\rm x,s})_{i,j} = \widetilde{\widetilde V}_{s,j,i} + \widetilde W_{*,s,j,i}
\end{eqnarray}
and
\begin{eqnarray}\label{Adef2dd}
&&(A_{\rm Har,\ell})_{i,j} + (A_{\rm x,\ell})_{i,j} \nonumber\\
&=&|\widetilde \theta|^2\widetilde{\widetilde V}_{\ell,j,i} + \widetilde\theta \sum_{m=1}^{N_H}c_{m,H}^*\widetilde V_{\ell,m,j,i} +\widetilde\theta^* \sum_{m=1}^{N_H}c_{m,H}\widetilde V^*_{\ell,m,j,i} + \sum_{m,n=1}^{N_H}c_{m,H}c^*_{n,H}V_{m,n,j,i} \nonumber\\
&& +   |\widetilde\theta|^2 \widetilde W_{*,\ell,j,i}  + \widetilde\theta^* \sum_{m=1}^{N_H} c_{m,H} \widetilde V^*_{\ell,j,m,i}  + \widetilde\theta \sum_{m=1}^{N_H} c_{m,H}^* \widetilde V_{\ell,i,j,m}
+ \sum_{m,n=1}^{N_H}c^*_{m,H}c_{n,H} V_{j,m,n,i}. \ \ \ \ \ \
\end{eqnarray}

Next, we begin to assemble the vector $b_{Hh}$.
The vector $b_{Hh}$ has the following expansion

\begin{eqnarray}\label{bhhdef}
(b_{Hh})_{i}
&=&\frac{1}{2}\int_\Omega\nabla\widetilde\phi_{\ell,h}(x)\nabla \psi_{i,H}(x)dx +\int_\Omega V_{\rm ext}\widetilde\phi_{\ell,h}(x)\psi_{i,H}(x)dx \nonumber\\
&&+ \int_\Omega V_{\rm Har}(\rho_{Hh})\widetilde\phi_{\ell,h}(x)\psi_{i,H}(x)dx +\int_{\Omega}(V_{\rm x}(P_{Hh})\widetilde\phi_{\ell,h})\psi_{i,H}(x)dx \nonumber\\
&=&\frac{1}{2}\int_\Omega\nabla\widetilde\phi_{\ell,h}(x)\nabla \psi_{i,H}(x)dx +\int_\Omega V_{\rm ext}\widetilde\phi_{\ell,h}(x)\psi_{i,H}(x)dx \nonumber\\
&&+ \sum_{s=1}^N\int_\Omega\int_\Omega \frac{ |\phi_{s,Hh}(y)|^2}{|x-y|}dy\widetilde\phi_{\ell,h}(x)\psi_{i,H}(x)dx \nonumber\\
&&+ \sum_{s=1}^N\int_\Omega\int_\Omega \frac{ \phi^*_{s,Hh}(y)\widetilde\phi_{\ell,h}(y)}{|x-y|}dy\phi_{s,Hh}(x)\psi_{i,H}(x)dx.
\end{eqnarray}

Denote
\begin{eqnarray}
(b_{\rm Kin,\ell})_{i}&=&\frac{1}{2}\int_\Omega\nabla\widetilde\phi_{\ell,h}(x)\nabla \psi_{i,H}(x)dx +\int_\Omega V_{\rm ext}\widetilde\phi_{\ell,h}(x)\psi_{i,H}(x)dx,  \label{bbkin}\\
(b_{\rm Har,s,\ell})_{i}&=&\int_\Omega\int_\Omega \frac{|\phi_{s,Hh}(y)|^2}{|x-y|}dy\widetilde\phi_{\ell,h}(x)\psi_{i,H}(x)dx, \label{bbhar}  \\
(b_{\rm x,s,\ell})_{i}&=&\int_\Omega\int_\Omega \frac{ \phi^*_{s,Hh}(y)\widetilde\phi_{\ell,h}(y)}{|x-y|}dy\phi_{s,Hh}(x)\psi_{i,H}(x)dx.  \label{bbx}
\end{eqnarray}
It is obvious that $b_{\rm Kin,\ell}$ remains unchanged during the correction step. For $b_{\rm Har,s,\ell}$ and $b_{\rm x,s,\ell}$, when $s\neq \ell$, we have
\begin{eqnarray}
&&(b_{\rm Har,s,\ell})_{i}+(b_{\rm x,s,\ell})_{i}\nonumber\\
&=&\int_\Omega\int_\Omega \frac{|\widetilde\phi_{s,h}(y)|^2}{|x-y|}dy\widetilde\phi_{\ell,h}(x)\psi_{i,H}(x)dx
 + \int_\Omega\Big\{\int_\Omega \frac{\widetilde\phi_{s,h}^*(y)\widetilde\phi_{\ell,h}(y) }{|x-y|}dy  \widetilde\phi_{s,h}(x)\psi_{i,H}(x)\Big\}dx \nonumber\\
&=&\int_\Omega \widetilde{\widetilde V}_{s,s}(x)  \widetilde\phi_{\ell,h}(x)\psi_{i,H}(x)dx
 + \int_\Omega  \widetilde{\widetilde V}_{\ell,s}(x)   \widetilde\phi_{s,h}(x)\psi_{i,H}(x)dx.
\end{eqnarray}
When $s=\ell$, we have
\begin{eqnarray}
&&(b_{\rm Har,\ell,\ell})_{i}+(b_{\rm x,\ell,\ell})_{i}\nonumber\\
&=&\int_\Omega\int_\Omega \frac{\Big|\widetilde \theta \widetilde\phi_{\ell,h}(y)+ \sum_{m=1}^{N_H}c_{m,H}\psi_{m,H}(y)\Big|^2}{|x-y|}dy\widetilde\phi_{\ell,h}(x)\psi_{i,H}(x)dx \nonumber\\
&& + \int_\Omega\Big\{\int_\Omega \frac{\Big(\widetilde \theta \widetilde\phi_{\ell,h}(y)+ \sum_{m=1}^{N_H}c_{m,H}\psi_{m,H}(y)\Big)^*\widetilde\phi_{\ell,h}(y)}{|x-y|}dy  \nonumber\\
&& \ \ \ \ \ \ \Big(\widetilde \theta\widetilde\phi_{\ell,h}(x)+ \sum_{m=1}^{N_H}c_{m,H}\psi_{m,H}(x)\Big)\psi_{i,H}(x)\Big\}dx \nonumber\\
&=&\int_\Omega \Big\{\Big(|\widetilde \theta|^2\widetilde{\widetilde V}_{\ell,\ell}(x) + \widetilde \theta\sum_{m=1}^{N_H} c_{m,H}^*\widetilde V_{\ell,m}(x)  + \widetilde \theta^*\sum_{m=1}^{N_H} c_{m,H}\widetilde V_{\ell,m}^*(x)
+ \sum_{m,n=1}^{N_H}c_{m,H}c_{n,H}^* V_{m,n}(x) \Big)  \nonumber\\
&& \ \ \ \ \ \ \widetilde\phi_{\ell,h}(x)\psi_{i,H}(x)\Big\}dx \nonumber\\
&& + \int_\Omega \Big(\widetilde \theta^* \widetilde{\widetilde V}_{\ell,\ell}(x) +\sum_{m=1}^{N_H} c^*_{m,H}\widetilde V_{\ell,m}(x)  \Big) \Big(\widetilde\theta\widetilde\phi_{\ell,h}(x)+ \sum_{m=1}^{N_H}c_{m,H}\psi_{m,H}(x)\Big)\psi_{i,H}(x)dx.
\end{eqnarray}

Define
\begin{eqnarray}\label{vttt}
\widetilde{\widetilde{\widetilde V}}_{s,m,\ell,i} = \int_\Omega\widetilde{\widetilde V}_{s,m}(x)\widetilde\phi_{\ell,h}(x)\psi_{i,H}(x)dx.
\end{eqnarray}
It is obvious that $\widetilde{\widetilde{\widetilde V}}_{s,m,\ell,i}$ remains unchanged during the correction step.
Then $b_{\rm Har,\ell}$ and $b_{\rm x,\ell}$ can be assembled in the following way:
\begin{eqnarray}\label{Adef2slbdg}
(b_{\rm Har,s,\ell})_{i}+(b_{\rm x,s,\ell})_{i}
= \widetilde{\widetilde{\widetilde V}}_{s,s,\ell,i}+\widetilde{\widetilde{\widetilde V}}_{\ell,s,s,i}
\end{eqnarray}
and
\begin{eqnarray}\label{Adef2bdg}
&&(b_{\rm Har,\ell,\ell})_{i}+(b_{\rm x,\ell,\ell})_{i}\nonumber\\
&=& |\widetilde \theta|^2\widetilde{\widetilde{\widetilde V}}_{\ell,\ell,\ell,i} + \widetilde\theta\sum_{m=1}^{N_H}c^*_{m,H}\widetilde W_{\ell,m,i}+ \widetilde\theta^*\sum_{m=1}^{N_H}c_{m,H}\widetilde W_{*,\ell,m,i} +\sum_{m,n=1}^{N_H}c_{m,H}c^*_{n,H}\widetilde V_{\ell,i,m,n} \nonumber\\
&& +  |\widetilde \theta|^2\widetilde{\widetilde{\widetilde V}}_{\ell,\ell,\ell,i}  + \widetilde\theta^* \sum_{m=1}^{N_H} c_{m,H} \widetilde{\widetilde V}_{\ell,m,i}
+ \widetilde\theta \sum_{m=1}^{N_H} c^*_{m,H} \widetilde W_{\ell,m,i} + \sum_{m,n}c^*_{m,H}c_{n,H} \widetilde V_{\ell,m,n,i}.
\end{eqnarray}

Finally, we come to assemble the scalar $\beta$.
The scalar $\beta$ has the following expansion
\begin{eqnarray}\label{Adefbeta}
\beta&=&\frac{1}{2}\int_\Omega\nabla\widetilde\phi_{\ell,h}(x)\nabla \widetilde\phi_{\ell,h}^*(x)dx +\int_\Omega V_{\rm ext}\widetilde\phi_{\ell,h}(x)\widetilde\phi_{\ell,h}^*(x)dx \nonumber\\
&&+ \int_\Omega V_{\rm Har}(\rho_{Hh})\widetilde\phi_{\ell,h}(x)\widetilde\phi_{\ell,h}^*(x)dx +\int_{\Omega}(V_{\rm x}(P_{Hh})\widetilde\phi_{\ell,h})\widetilde\phi_{\ell,h}^*(x)dx \nonumber\\
&=&\frac{1}{2}\int_\Omega\nabla\widetilde\phi_{\ell,h}(x)\nabla \widetilde\phi_{\ell,h}^*(x)dx +\int_\Omega V_{\rm ext}\widetilde\phi_{\ell,h}(x)\widetilde\phi_{\ell,h}^*(x)dx \nonumber\\
&&+ \sum_{s=1}^N\int_\Omega\int_\Omega \frac{ |\phi_{s,Hh}(y)|^2}{|x-y|}dy\widetilde\phi_{\ell,h}(x)\widetilde\phi_{\ell,h}^*(x)dx \nonumber\\
&&+ \sum_{s=1}^N\int_\Omega\int_\Omega \frac{ \phi^*_{s,Hh}(y)\widetilde\phi_{\ell,h}(y)}{|x-y|}dy\phi_{s,Hh}(x)\widetilde\phi_{\ell,h}^*(x)dx.\ \
\end{eqnarray}

Denote
\begin{eqnarray}
\beta_{\rm Kin,\ell}&=&\frac{1}{2}\int_\Omega\nabla\widetilde\phi_{\ell,h}(x)\nabla \widetilde\phi_{\ell,h}^*(x)dx +\int_\Omega V_{\rm ext}\widetilde\phi_{\ell,h}(x)\widetilde\phi_{\ell,h}^*(x)dx, \label{betakin}\\
\beta_{\rm Har,s,\ell}&=&\int_\Omega\int_\Omega \frac{|\phi_{s,Hh}(y)|^2}{|x-y|}dy\widetilde\phi_{\ell,h}(x)\widetilde\phi_{\ell,h}^*(x)dx,\\
\beta_{\rm x,s,\ell}&=&\int_\Omega\int_\Omega \frac{ \phi^*_{s,Hh}(y)\widetilde\phi_{\ell,h}(y)}{|x-y|}dy\phi_{s,Hh}(x)\widetilde\phi_{\ell,h}^*(x)dx.
\end{eqnarray}
It is obvious that $\beta_{\rm Kin,\ell}$ remains unchanged during the correction step. For $\beta_{\rm Har,s,\ell}$ and $\beta_{\rm x,s,\ell}$, when $s\neq \ell$, we have
\begin{eqnarray}
&&\beta_{\rm Har,s,\ell}+\beta_{\rm x,s,\ell}\nonumber\\
&=&\int_\Omega\int_\Omega \frac{|\widetilde\phi_{s,h}(y)|^2}{|x-y|}dy\widetilde\phi_{\ell,h}(x)\widetilde\phi_{\ell,h}^*(x)dx
 + \int_\Omega\int_\Omega \frac{\widetilde\phi_{s,h}^*(y)\widetilde\phi_{\ell,h}(y)}{|x-y|}dy
\widetilde\phi_{s,h}(x)\widetilde\phi_{\ell,h}^*(x)dx \nonumber\\
&=&\int_\Omega  \widetilde{\widetilde V}_{s,s}(x) \widetilde\phi_{\ell,h}(x)\widetilde\phi_{\ell,h}^*(x)dx
 + \int_\Omega \widetilde{\widetilde V}_{\ell,s}(x) \widetilde\phi_{s,h}(x)\widetilde\phi_{\ell,h}^*(x)dx
\end{eqnarray}
and
\begin{eqnarray}
&&\beta_{\rm Har,\ell,\ell}+\beta_{\rm x,\ell,\ell}\nonumber\\
&=&\int_\Omega\int_\Omega \frac{\Big|\widetilde \theta \widetilde\phi_{\ell,h}(y)+ \sum_{m=1}^{N_H}c_{m,H}\psi_{m,H}(y)\Big|^2}{|x-y|}dy\widetilde\phi_{\ell,h}(x)\widetilde\phi_{\ell,h}^*(x)dx \nonumber\\
&& + \int_\Omega\Big\{\int_\Omega \frac{\Big(\widetilde \theta \widetilde\phi_{\ell,h}(y)+ \sum_{m=1}^{N_H}c_{m,H}\psi_{m,H}(y)\Big)^*\widetilde\phi_{\ell,h}(y)}{|x-y|}dy \nonumber\\
&&\ \ \ \ \ \ \ \ \ \ \ \ \Big(\widetilde \theta\widetilde\phi_{\ell,h}(x)+ \sum_{m=1}^{N_H}c_{m,H}\psi_{m,H}(x)\Big)\widetilde\phi_{\ell,h}^*(x)\Big\}dx \nonumber\\
&=&\int_\Omega \Big\{ \Big(|\widetilde \theta|^2\widetilde{\widetilde V}_{\ell,\ell}(x) + \widetilde \theta\sum_{m=1}^{N_H} c_{m,H}^*\widetilde V_{\ell,m}(x)  + \widetilde \theta^*\sum_{m=1}^{N_H} c_{m,H}\widetilde V_{\ell,m}^*(x)
+ \sum_{m,n=1}^{N_H}c_{m,H}c_{n,H}^* V_{m,n}(x) \Big) \nonumber\\
&& \ \ \ \ \ \ \ \widetilde\phi_{\ell,h}(x)\widetilde\phi_{\ell,h}^*(x)\Big\}dx \nonumber\\
&& + \int_\Omega \Big(\widetilde \theta^* \widetilde{\widetilde V}_{\ell,\ell}(x) +\sum_{m=1}^{N_H} c_{m,H}^*\widetilde V_{\ell,m}(x)  \Big) \Big(\widetilde\theta\widetilde\phi_{\ell,h}(x)+ \sum_{m=1}^{N_H}c_{m,H}\psi_{m,H}(x)\Big)\widetilde\phi_{\ell,h}^*(x)dx.
\end{eqnarray}

Define
\begin{eqnarray}\label{vtttt}
\widetilde{\widetilde{\widetilde{\widetilde V}}}_{s,m,j,\ell} = \int_\Omega\widetilde{\widetilde V}_{s,m}(x)\widetilde\phi_{j,h}(x)\widetilde\phi_{\ell,h}^*(x)dx.
\end{eqnarray}
It is obvious that  $\widetilde{\widetilde{\widetilde{\widetilde V}}}_{s,m,j,\ell}$ remains unchanged during the correction step.
Then $\beta$ can be assembled in the following way
\begin{eqnarray}\label{Adef2sjbet}
\beta_{\rm Har,s,\ell}+\beta_{\rm x,s,\ell} = \widetilde{\widetilde{\widetilde{\widetilde V}}}_{s,s,\ell,\ell} +\widetilde{\widetilde{\widetilde{\widetilde V}}}_{\ell,s,s,\ell}
\end{eqnarray}
and
\begin{eqnarray}\label{Adef2bet}
&&\beta_{\rm Har,\ell,\ell}+\beta_{\rm x,\ell,\ell}\nonumber\\
&=& |\widetilde \theta|^2\widetilde{\widetilde{\widetilde{\widetilde V}}}_{\ell,\ell,\ell,\ell} + \widetilde\theta\sum_{m=1}^{N_H}c_{m,H}^*\widetilde{\widetilde{\widetilde V}}_{\ell,\ell,\ell,m}
+ \widetilde\theta^*\sum_{m=1}^{N_H}c_{m,H}\widetilde{\widetilde{\widetilde V}}_{\ell,\ell,\ell,m}^*
+ \sum_{m,n=1}^{N_H}c_{m,H}c_{n,H}^*\widetilde{\widetilde V}_{\ell,m,n} \nonumber\\
&& +  |\widetilde \theta|^2\widetilde{\widetilde{\widetilde{\widetilde V}}}_{\ell,\ell,\ell,\ell}  + \widetilde\theta^* \sum_{m=1}^{N_H} c_{m,H} \widetilde{\widetilde{\widetilde V}}_{\ell,\ell,\ell,m}^*
+ \widetilde\theta \sum_{m=1}^{N_H} c_{m,H}^* \widetilde{\widetilde{\widetilde V}}_{\ell,\ell,\ell,m}
+  \sum_{m,n=1}^{N_H}c_{m,H}c_{n,H}^* \widetilde W_{*,\ell,m,n}. \ \ \ \ \ \ \
\end{eqnarray}
As we can see, to construct the stiffness matrix in (\ref{Eigenvalue_H_h}), we just need to compute:
\begin{enumerate}
\item Two tensors $\mathbf V$ and $\mathbf{\widetilde V}_\ell$ with $(\mathbf V)_{m,n,i,j}=V_{m,n,i,j}$ and $(\mathbf{\widetilde V}_\ell)_{m,i,j}=\widetilde V_{\ell,m,i,j}$.
\item Three matrices $\mathbf{\widetilde{\widetilde V}_{\ell}}$, $\mathbf {\widetilde W}_\ell$, $\mathbf {\widetilde W}_{*,\ell}$ with $(\mathbf{\widetilde{\widetilde V}_{\ell}})_{i,j}=\widetilde{\widetilde V}_{\ell,i,j}$, $(\mathbf {\widetilde W}_\ell)_{i,j} = \widetilde W_{\ell,i,j}$, $(\mathbf {\widetilde W}_{*,\ell})_{i,j} = \widetilde W_{*,\ell,i,j}$.
\item One vector $\mathbf{\widetilde{\widetilde{\widetilde V}}}_{\ell}$ with $(\mathbf{\widetilde{\widetilde{\widetilde V}}}_{\ell})_{i}=\widetilde{\widetilde{\widetilde V}}_{\ell,\ell,\ell,i}$.
\item  One scalar $\widetilde{\widetilde{\widetilde{\widetilde V}}}_{\ell}$ with $\widetilde{\widetilde{\widetilde{\widetilde V}}}_{\ell}=\widetilde{\widetilde{\widetilde{\widetilde V}}}_{\ell,\ell,\ell,\ell}$.
\item[]\hspace*{-\leftmargini} Besides, we need to compute some matrices and vectors that depend on the other wavefunctions:
\item Two matrices $\mathbf{\widetilde{\widetilde{\widetilde V}}_{1\ell}}$, $\mathbf{\widetilde{\widetilde{\widetilde V}}_{2\ell}}$
with $(\mathbf{\widetilde{\widetilde{\widetilde V}}_{1\ell}})_{s,i}=\widetilde{\widetilde{\widetilde V}}_{s,s,\ell,i}$, $(\mathbf{\widetilde{\widetilde{\widetilde V}}_{2\ell}})_{s,i}=\widetilde{\widetilde{\widetilde V}}_{\ell,s,s,i}$.
\item Two vectors $\mathbf{\widetilde{\widetilde{\widetilde{\widetilde V}}}}_{1\ell}$ and $\mathbf{\widetilde{\widetilde{\widetilde{\widetilde V}}}}_{2\ell}$
with $(\mathbf{\widetilde{\widetilde{\widetilde{\widetilde V}}}}_{1\ell})_{s}=\widetilde{\widetilde{\widetilde{\widetilde V}}}_{s,s,\ell,\ell}$
and $(\mathbf{\widetilde{\widetilde{\widetilde{\widetilde V}}}}_{2\ell})_{s}=\widetilde{\widetilde{\widetilde{\widetilde V}}}_{\ell,s,s,\ell}$.
\end{enumerate}

For these quantities, the tensor $\mathbf V$ remains unchanged throughout the algorithm,  thus $\mathbf V$ only needs to be computed one time at the beginning of the algorithm.
The other quantities remain unchanged in the SCF iteration of each correction space, thus these quantities only need to be computed one time before the SCF iteration in each correction space.

Based on the definitions of these quantities, the main point is to compute the involved functions $\widetilde{\widetilde V}_{\ell,s}(x)$, $\widetilde V_{\ell,m}(x)$ and $V_{m,n}(x)$.

To compute the functions $\widetilde{\widetilde V}_{\ell,s}(x)$, $\widetilde V_{\ell,m}(x)$ and $V_{m,n}(x)$, we need to solve some Poisson equations as follows:
\begin{eqnarray}\label{type1}
-\Delta \widetilde{\widetilde V}_{\ell,s}(x) = 4\pi \widetilde\phi_{\ell,h}\widetilde\phi_{s,h}^*, \ \ \ s=1,\cdots,N,
\end{eqnarray}
\begin{eqnarray}\label{type2}
-\Delta \widetilde V_{\ell,m}(x) = 4\pi \widetilde\phi_{\ell,h}\psi_{m,H}, \ \ \ m=1,\cdots,N_H,
\end{eqnarray}
\begin{eqnarray}\label{type3}
-\Delta V_{m,n}(x) = 4\pi \psi_{m,H}\psi_{n,H}, \ \ \ m,n=1,\cdots,N_H.
\end{eqnarray}

Because the coarse space $V_H$ is fixed, the equations (\ref{type3}) only need to be solved one time for the entire algorithm.
These functions can even be stored in the package in advance.
The equations (\ref{type1}) and (\ref{type2}) only need to be solved before the SCF iteration in each correction space.

\begin{remark}

In the multilevel correction type of AFEM, we need to solve some linearized boundary value problems in adaptive finite element spaces.
For these boundary value problems,
the exchange potentials are linearized by the approximate wavefunctions derived from the previous mesh level.
Then we will obtain a sparse system of equations that can be solved efficiently. Next, we will briefly describe how to assemble the vector corresponding to the exchange potential.
In fact, after deriving the approximate wavefunctions
based on the above mentioned process, we can assemble $V_{\rm x}(P_{h_{k+1}})$ on $V_{H,h_{k+1}}$ first, and then prolong the exchange potential from $V_{H,h_{k+1}}$  to $V_{h_{k+2}}$.
This will form the righthand side term of the linearized boundary value problem defined in $V_{h_{k+2}}$.

Assume the final approximate solution derived in $V_{H,h_{k+1}}$ is denoted by
\begin{eqnarray*}
\phi_{\ell,Hh}=\widetilde \theta_\ell \widetilde\phi_{\ell,h}+ \sum_{m=1}^{N_H}c_{m,H,\ell}\psi_{m,H}.
\end{eqnarray*}
To discretize the exchange potential, the key point is to derive $\int_\Omega\frac{\phi_{\ell,Hh}(y)\phi_{s,Hh}(y)}{|x-y|}dy$. By simple calculations, we can derive
\begin{eqnarray*}
&&\int_\Omega\frac{\phi_{\ell,Hh}(y)\phi_{s,Hh}^*(y)}{|x-y|}dy \nonumber\\
&=& \int_\Omega \frac{\big(\widetilde \theta_\ell \widetilde\phi_{\ell,h}(y) +\sum_{m=1}^{N_H} c_{m,H,\ell}\psi_{m,H}(y)  \big) \big(\widetilde\theta_s^*\widetilde\phi_{s,h}^*(y)
+ \sum_{m=1}^{N_H}c_{m,H,s}^*\psi_{m,H}(y)\big)}{|x-y|}dy \nonumber\\
&=& \widetilde \theta_\ell \widetilde \theta_s^* \widetilde{\widetilde V}_{\ell,s} + \sum_{m=1}^{N_H} c_{m,H,\ell}\widetilde\theta_s^*\widetilde V_{s,m}^*+ \sum_{m=1}^{N_H} c_{m,H,s}^*\widetilde\theta_\ell\widetilde V_{\ell,m}+\sum_{m,n=1}^{N_H}c_{m,H,\ell}c_{n,H,s}^*V_{m,n}.
\end{eqnarray*}
We can see that all the required components have already been derived when we solve the small-scale Hartree--Fock equation in the correction space. Thus, only a few matrix and tensor operations are required to discretize the exchange potential.
\end{remark}

In Algorithm \ref{multilevelAFEMim}, based on above discussion and preparation, we propose the efficient
implementation strategy for Algorithm \ref{multilevelAFEM}. 

\begin{algorithm}\label{multilevelAFEMim}Efficient implementation strategy for Algorithm \ref{multilevelAFEM}
\begin{enumerate}
\item Compute the matrices $M_H$ as (\ref{massmatrix}), $A_{\rm Kin}$ as (\ref{akin}), compute $V_{m,n}$ according to (\ref{type3}) and then compute the tensor $\mathbf{V}$ as (\ref{vmnjk}) in the coarse space $V_H$.

\item
Solve the Hartree--Fock equation (\ref{f1Hartree--Fock}) in the initial finite element space $V_{h_1}$.
\item
Set $k=1$.
\item
Compute the local error indicators ${\eta}_{h_k}(\{\lambda_{\ell,h_k},\phi_{\ell,h_k}\}_{\ell=1}^N,T)$ for each element $T\in\mathcal T_{h_k}$.
\item
Construct $\mathcal M_{h_k}\subset \mathcal T_{h_k}$ by D\"{o}rfler's marking strategy and
refine $\mathcal T_{h_k}$ to get a new mesh $\mathcal T_{h_{k+1}}$.
\item Solve the linear boundary value problem (\ref{Aux_Linear_Problem}) to derive $\widetilde{\phi}_{\ell,h_{k+1}} \in V_{h_{k+1}},\ \ell=1,\cdots,N$.
\item For $\ell=1,\cdots,N$:
\begin{enumerate}[label=(\Alph*)]
\item Define $V_{H,h_{k+1}}=V_H+{\rm span}\{\widetilde{\phi}_{\ell,h_{k+1}}\}$ and  compute the following quantities:
\begin{enumerate}[label=(\arabic*)]
\item Compute $\widetilde{\widetilde V}_{\ell,s}(x), s=1,\cdots,N$ and $\widetilde V_{\ell,m}(x), m=1,\cdots,N_H$ as (\ref{type1}) and (\ref{type2}).
\item Compute the tensor $\mathbf{\widetilde V}_\ell$ as (\ref{lmkj}), three matrices $\mathbf{\widetilde{\widetilde V}_{\ell}}$, $\mathbf {\widetilde W}_\ell$,
$\mathbf {\widetilde W}_{*,\ell}$ as (\ref{vlkj})  (\ref{wlkj}) and (\ref{*wlkj}),
a vector $\mathbf{\widetilde{\widetilde{\widetilde V}}}_{\ell}$ as (\ref{vttt}) and a scalar $\widetilde{\widetilde{\widetilde{\widetilde V}}}_{\ell}$ as (\ref{vtttt}) .
\item Compute two matrices $\mathbf{\widetilde{\widetilde{\widetilde V}}_{1\ell}}$, $\mathbf{\widetilde{\widetilde{\widetilde V}}_{2\ell}}$
as (\ref{vttt}), two vectors $\mathbf{\widetilde{\widetilde{\widetilde{\widetilde V}}}}_{1\ell}$ and $\mathbf{\widetilde{\widetilde{\widetilde{\widetilde V}}}}_{2\ell}$
as (\ref{vtttt}).
\item Compute the vector $b_{\rm Kin}$ as (\ref{bbkin}) and the scalar $\beta_{\rm Kin}$ as (\ref{betakin}).
\item Compute the vector $c_{Hh}$ and the scalar $\gamma$ as (\ref{massmatrix}).
\end{enumerate}
\item Do the following loop:
\begin{enumerate}[label=(\arabic*)]
\item Compute $A_{\rm Har,s}+A_{\rm x,s}$ as (\ref{Adef2sldd}) and  (\ref{Adef2dd}). Then $A_H=A_{\rm Kin}+\sum_s(A_{\rm Har,s}+A_{\rm x,s})$.
\item Compute $b_{\rm Har,s,\ell}+b_{\rm x,s,\ell}$ as (\ref{Adef2slbdg}) and (\ref{Adef2bdg}). Then $b_{Hh}=b_{\rm Kin}+\sum_s(b_{\rm Har,s,\ell}+b_{\rm x,s,\ell})$.
\item Compute $\beta_{\rm Har,s,\ell}+\beta_{\rm x,s,\ell}$ as (\ref{Adef2sjbet}) and (\ref{Adef2bet}). Then $\beta=\beta_{\rm Kin}+\sum_s(\beta_{\rm Har,s,\ell}+\beta_{\rm x,s,\ell})$.
\item Solve the eigenvalue problem (\ref{Eigenvalue_H_h}).
\item If the desired accuracy is satisfied, stop. Else, go to step (1) of (B) and continue.
\end{enumerate}
\end{enumerate}
\item
Let $k=k+1$ and go to step 4.
\end{enumerate}
\end{algorithm}

\begin{remark}
As we can see, during the SCF iteration in the $7 (B)$ step of Algoirthm \ref{multilevelAFEMim}, we just need to perform some small-scale matrix and tensor calculations. No additional matrix assembling are needed.
Thus, the correction step can be implemented efficiently.

Next, we briefly analyze the computational work for solving the $\ell$-th wavefunction.
Assume that the linearized boundary value problems can be solved with linear computational work based on some efficient numerical methods such as the multigrid method.
Then for the boundary value problem (\ref{Aux_Linear_Problem}), the computational work is $O(N_k)$.

For the boundary value problems (\ref{type1}), the computational work is $O(NN_k)$.

For the boundary value problems (\ref{type2}), the computational work is $O(N_HN_k)$.

For the  tensor $\mathbf{\widetilde V}_\ell$ as (\ref{lmkj}), the computational work is $O(N_HN_k)$.

For the matrix $\mathbf{\widetilde{\widetilde V}_{\ell}}$ as (\ref{vlkj}), the computational work is $O(N_k)$.

For the  matrices $\mathbf {\widetilde W}_\ell$,
$\mathbf {\widetilde W}_{*,\ell}$ as   (\ref{wlkj}) and (\ref{*wlkj}), the computational work is $O(N_HN_k)$.

For the vector $\mathbf{\widetilde{\widetilde{\widetilde V}}}_{\ell}$ as (\ref{vttt}) and a scalar $\widetilde{\widetilde{\widetilde{\widetilde V}}}_{\ell}$ as (\ref{vtttt}), the computational work is $O(N_k)$.

For the two matrices $\mathbf{\widetilde{\widetilde{\widetilde V}}_{1\ell}}$, $\mathbf{\widetilde{\widetilde{\widetilde V}}_{2\ell}}$
as (\ref{vttt}), the computational work is $O(NN_k)$.

For the two vectors $\mathbf{\widetilde{\widetilde{\widetilde{\widetilde V}}}}_{1\ell}$ and $\mathbf{\widetilde{\widetilde{\widetilde{\widetilde V}}}}_{2\ell}$
as (\ref{vtttt}), the computational work is $O(NN_k)$.

For the vector $b_{\rm Kin}$ as (\ref{bbkin}) and the scalar $\beta_{\rm Kin}$ as (\ref{betakin}), the computational work is $O(N_k)$.

For the vector $c_{Hh}$ and the scalar $\gamma$ as (\ref{massmatrix}), the computational work is $O(N_k)$.

During the SCF iteration, we need to perform some matrix operations:

 For $A_H=A_{\rm Kin}+\sum_s(A_{\rm Har,s}+A_{\rm x,s})$,  the computational work is $O(NN_H^2+N_H^3)$.

 For  $b_{Hh}=b_{\rm Kin}+\sum_s(b_{\rm Har,s,\ell}+b_{\rm x,s,\ell})$,
 the computational work is $O(NN_H+N_H^3)$.

 For  $\beta=\beta_{\rm Kin}+\sum_s(\beta_{\rm Har,s,\ell}+\beta_{\rm x,s,\ell})$,
the computational work is $O(N+N_H^2)$.

Assume $\omega$ times SCF iteration are needed and the computational work for small-scale linear eigenvalue problem (\ref{Eigenvalue_H_h}) is $O(M_H)$. Then the total computational work is
\begin{eqnarray}\label{comest}
O\big((N+N_H)N_k+\omega(NN_H^2+N_H^3+M_H)\big).
\end{eqnarray}
For a given molecular system, to derive a more accurate approximation, the value of $N_k$ is going to get increasingly larger. However, the number of electrons and the dimension of correction space remain fixed.
In other words, the values of $N$, $N_H$ and $M_H$ remain unchanged.
In addition, the coefficient of $N_k$ in (\ref{comest}) is independent from the iteration time $\omega$.
Thus, a linear computational work can be obtained with the refinement of mesh, which is gradually freed from the influence of the number of SCF iterations.

\end{remark}

\begin{remark}
In terms of storage, direct AFEM needs to store large-scale dense matrices which involves $O(N_k^2)$ data.
In the proposed algorithm, the involved dense matrices are $\mathbf {\widetilde W}_\ell$,
$\mathbf {\widetilde W}_{*,\ell}$, $\mathbf{\widetilde{\widetilde{\widetilde V}}_{1\ell}}$, $\mathbf{\widetilde{\widetilde{\widetilde V}}_{2\ell}}$.
The dimension dim($\mathbf {\widetilde W}_\ell$)=
dim($\mathbf {\widetilde W}_{*,\ell}$)=$N_H\times N_H$, dim($\mathbf{\widetilde{\widetilde{\widetilde V}}_{1\ell}}$)=dim($\mathbf{\widetilde{\widetilde{\widetilde V}}_{2\ell}}$)$=N\times N_H$.
Thus, the amount of data to be stored are $O(NN_H+N_H^2)$. The total number of stored data are $O(N_k+NN_H+N_H^2)$, which is a significant improvement over the direct AFEM.
\end{remark}

\begin{remark}
Finally, we highlight that Algorithm \ref{multilevelAFEMim} is inherently suitable for eigenpairwise parallel computing due to its distinctive structure.
Notably, Algorithm \ref{multilevelAFEMim} handles different orbits independently, thus circumventing the need for large-scale inner product computation
for orthogonalization, which often poses a bottleneck in parallel computing.
Consequently, this strategy facilitates parallel computing by allowing different wavefunctions to be processed concurrently on separate processors.
\end{remark}

\section{Numerical experiments}
In this section, we present several numerical examples to validate the
efficiency of the proposed numerical method.
These numerical examples are performed on a cluster with 90 nodes.
Each computing node has two 20-core with Intel Xeon E5-2660 v3 processors at 2.6 GHz and 192 GB memory.

In our numerical experiments, we consider the Hartree--Fock equation for Lithium (HLi), Methane (CH$_4$), Benzene(C$_6$H$_6$) and Ethanol(C$_2$H$_6$O).
To exhibit the advantages of the proposed algorithm, we compare the numerical results with
the direct adaptive finite element method (i.e. solving Hartree--Fock equation by Algorithm \ref{Multilevel Correction1} directly in each adaptive finite element space)
from four perspectives including precision, solving efficiency, memory consumption and parallel scalability.
For both Algorithm \ref{multilevelAFEMim} and direct AFEM, the linear finite element basis are employed.

\subsection{Precision}
In this subsection, we aim to exhibit the accuracy of Algorithm \ref{multilevelAFEMim}.
For this end, we present the approximate energies derived from Algorithm \ref{multilevelAFEMim} and the results based on Gaussian basis functions employing
the computational chemistry package NWChem \cite{valiev}  in Table \ref{timetableex1}.

From Table \ref{timetableex1}, we can find that Algorithm \ref{multilevelAFEMim} is able to derive quite accurate approximate solutions as that of  NWChem.

\begin{table}[htbp]
\begin{center}
\resizebox{0.6\textwidth}{!}{
\begin{tabular}{c | c c c }\hline
Atom & Energy of Algorithm \ref{multilevelAFEMim} & Energy of NWChem \\ \hline
Lithium hydride  & --7.9842  & --7.9842 \\ 
Methane          &  --40.1998 & --40.1996\\ 
Ethanol   &  --154.1057 & --154.1065\\ 
Benzene         & --230.7265 & --230.7284\\ \hline
\end{tabular}
}
\end{center}
\caption{The approximate energies derived by Algorithm \ref{multilevelAFEMim} and NWChem.}
\label{timetableex1}
\end{table}

\subsection{Solving efficiency}
In this subsection, we present the computational time (in seconds) of Algorithm \ref{multilevelAFEMim} and direct AFEM with the refinement of mesh.
The corresponding results and the acceleration ratio between the two algorithms are presented in Tables \ref{timetableex21}--\ref{timetableex24}.
The acceleration ratio is obtained by dividing the time of direct AFEM by the time of Algorithm \ref{multilevelAFEMim}.
The symbol ``-" means the computer runs out of memory.
From Tables \ref{timetableex21}--\ref{timetableex24}, we can find that Algorithm  \ref{multilevelAFEMim} has a significant advantage over the direct AFEM.
First, for the same molecule, we can see that the advantage increases with the refinement of mesh. Next, for different molecules,
we can see that the advantage becomes more apparent for more complex models.


%

\begin{table}[htbp]
\begin{center}
\resizebox{0.69\textwidth}{!}{
\begin{tabular}{c | c c c}\hline
Atom & Time of direct AFEM
& Time of Algorithm \ref{multilevelAFEMim} & Ratio \\ \hline
Lithium hydride  &1.4466E+5    & 1.6407E+2  & 8.8171E+2 \\ 
Methane       &  5.3731E+5     &  3.0857E+2 & 1.7412E+3\\ 
Ethanol & -   &  1.1323E+3 & -\\ 
Benzene       &  -     & 2.4201E+3 & -\\ \hline
\end{tabular}
}
\end{center}
\caption{The computational time (in seconds) of the direct AFEM, Algorithm \ref{multilevelAFEMim} with the same energy accuracy 1E-1 and the acceleration ratio.}
\label{timetableex21}
\end{table}

\begin{table}[htbp]
\begin{center}
\resizebox{0.69\textwidth}{!}{
\begin{tabular}{c | c c c}\hline
Atom & Time of direct AFEM
& Time of Algorithm \ref{multilevelAFEMim} & Ratio \\ \hline
Lithium hydride  &2.7947E+6    & 3.0526E+2  & 9.1551E+3 \\ 
Methane       &  1.0986E+7     & 5.8011E+2 & 1.8939E+4\\ 
Ethanol & -  &  2.1288E+3 & -\\ 
Benzene       &  -     & 4.5493E+3 & -\\ \hline
\end{tabular}
}
\end{center}
\caption{The computational time (in seconds) of the direct AFEM, Algorithm \ref{multilevelAFEMim} with the same energy accuracy 1E-2 and the acceleration ratio.}
\label{timetableex22}
\end{table}

\begin{table}[htbp]
\begin{center}
\resizebox{0.69\textwidth}{!}{
\begin{tabular}{c | c c c}\hline
Atom & Time of direct AFEM
& Time of Algorithm \ref{multilevelAFEMim} & Ratio \\ \hline
Lithium hydride  & -    & 1.1323E+3  & - \\ 
Methane       &  -     &  2.2079E+3 & -\\ 
Ethanol & -   &  7.8043E+3 & -\\ 
Benzene       &  -     & 1.6623E+4 & -\\ \hline
\end{tabular}
}
\end{center}
\caption{The computational time (in seconds) of the direct AFEM, Algorithm \ref{multilevelAFEMim} with the same energy accuracy 1E-3 and the acceleration ratio.}
\label{timetableex23}
\end{table}

\begin{table}[htbp]
\begin{center}
\resizebox{0.69\textwidth}{!}{
\begin{tabular}{c | c c c}\hline
Atom & Time of direct AFEM
& Time of Algorithm \ref{multilevelAFEMim} & Ratio \\ \hline
Lithium hydride  & -    & 2.4202E+3  & - \\ 
Methane       &  -     &  4.5560E+3 & - \\ 
Ethanol & -   &  1.6765E+4 & -\\ 
Benzene       &  -    & 3.5623E+4 & -\\ \hline
\end{tabular}
}
\end{center}
\caption{The computational time (in seconds) of the direct AFEM, Algorithm \ref{multilevelAFEMim} with the same energy accuracy 1E-4 and the acceleration ratio.}
\label{timetableex24}
\end{table}

\subsection{Memory consumption}
A major challenge in solving the Hartree--Fock equation is that traditional algorithms generate dense matrices, which consume a large amount of memory.
Our algorithm, however, avoids generating large-scale dense matrices,
resulting in significant memory savings. This result is verified by the numerical data presented in Tables \ref{timetableex31}--\ref{timetableex34},
which show that the memory requirement of our algorithm is greatly reduced compared to that of the direct AFEM.
Meanwhile, the memory usage analysis reveals that the algorithm's advantages become increasingly prominent when handling larger-scale computations and more complex models.

\begin{table}[htbp]
\begin{center}
\resizebox{0.75\textwidth}{!}{
\begin{tabular}{c | c c c}\hline
Atom & Memory of direct AFEM
& Memory of Algorithm \ref{multilevelAFEMim} & Ratio \\ \hline
Lithium hydride  &1.1601E+1    & 1.0767E$-$1  & 1.0774E+2 \\ 
Methane       &  7.6577E+1     & 1.4260E$-$1 & 5.3697E+2\\ 
Ethanol & -   &  1.5974E$-$1 & -\\ 
Benzene       & - & 2.3050E$-$1 & -\\ \hline
\end{tabular}
}
\end{center}
\caption{The memory (in GB) of the direct AFEM, Algorithm \ref{multilevelAFEMim} with the same energy accuracy 1E-1 and the consumption ratio.}
\label{timetableex31}
\end{table}

\begin{table}[htbp]
\begin{center}
\resizebox{0.75\textwidth}{!}{
\begin{tabular}{c | c c c}\hline
Atom & Memory of direct AFEM
& Memory of Algorithm \ref{multilevelAFEMim} & Ratio \\ \hline
Lithium hydride  &2.3001E+1   & 1.4887E$-$1  & 1.5450E+2 \\ 
Methane       &  1.8694E+2     & 1.7475E$-$1 & 1.0697E+3\\ 
Ethanol & -  &  2.1387E$-$1 & -\\ 
Benzene       &  -     & 2.4281E$-$1 & -\\ \hline
\end{tabular}
}
\end{center}
\caption{The memory (in GB) of the direct AFEM, Algorithm \ref{multilevelAFEMim} with the same energy accuracy 1E-2 and the consumption ratio.}
\label{timetableex32}
\end{table}

\begin{table}[htbp]
\begin{center}
\resizebox{0.75\textwidth}{!}{
\begin{tabular}{c | c c c}\hline
Atom & Memory of direct AFEM
& Memory of Algorithm \ref{multilevelAFEMim} & Ratio \\ \hline
Lithium hydride  & -    & 2.6732E$-$1  & - \\ 
Methane       &  -     &  3.2284E$-$1 & -\\ 
Ethanol & -   &  4.0338E$-$1 & -\\ 
Benzene       &  -     & 5.6582E$-$1 & -\\ \hline
\end{tabular}
}
\end{center}
\caption{The memory (in GB) of the direct AFEM, Algorithm \ref{multilevelAFEMim} with the same energy accuracy 1E-3 and the consumption ratio.}
\label{timetableex33}
\end{table}

\begin{table}[htbp]
\begin{center}
\resizebox{0.75\textwidth}{!}{
\begin{tabular}{c | c c c}\hline
Atom & Memory of direct AFEM
& Memory of Algorithm \ref{multilevelAFEMim} & Ratio \\ \hline
Lithium hydride  & -    & 8.9169E$-$1  & - \\ 
Methane       &  -     &  1.1426E+1 & -\\ 
Ethanol & -   &  1.3787E+1 & -\\ 
Benzene       &  -    & 2.0454E+1 & -\\ \hline
\end{tabular}
}
\end{center}
\caption{The memory (in GB) of the direct AFEM, Algorithm \ref{multilevelAFEMim} with the same energy accuracy 1E-4 and the consumption ratio.}
\label{timetableex34}
\end{table}

\subsection{Parallel scalability}
\begin{figure}[htbp]
\centering
\includegraphics[width=5.0cm]{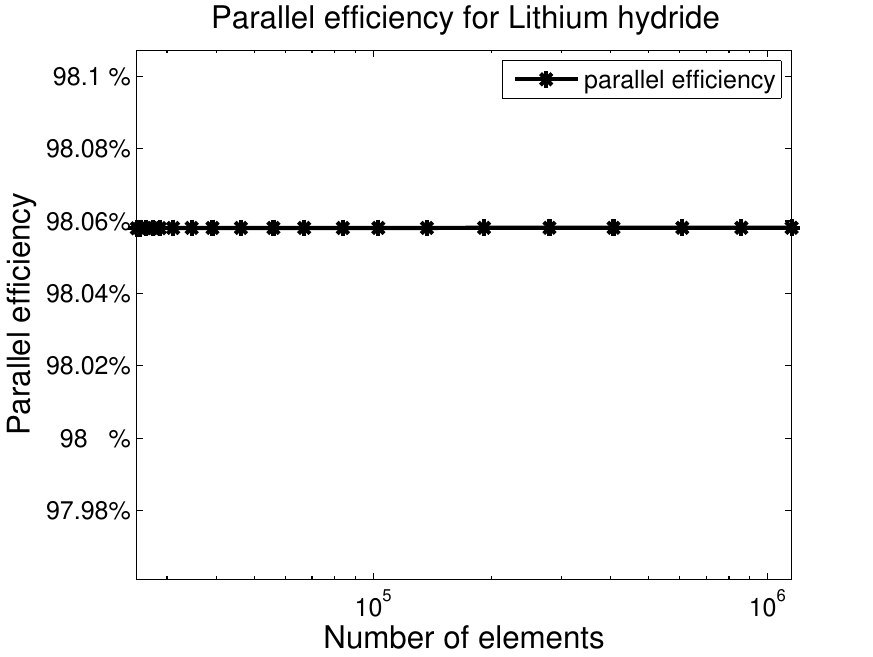}\ \ \ \ \ \ \ \ \
\includegraphics[width=5.0cm]{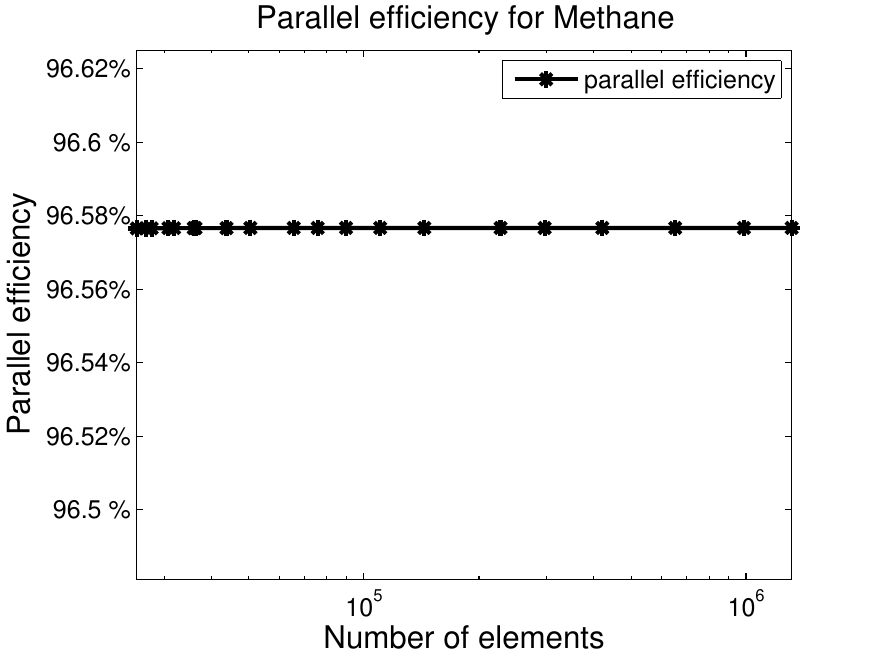}\\
\includegraphics[width=5.0cm]{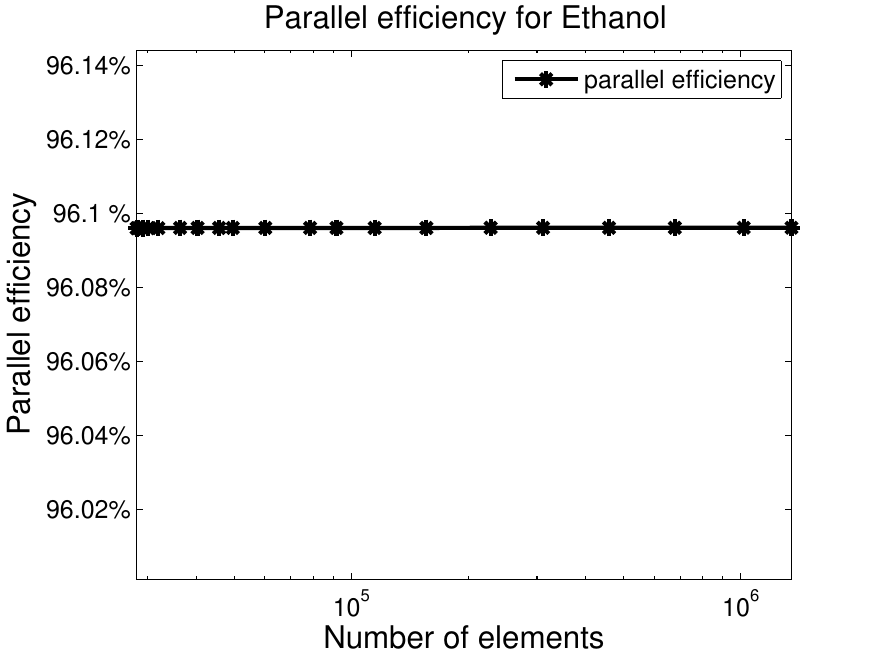}\ \ \ \ \ \ \ \ \
\includegraphics[width=5.0cm]{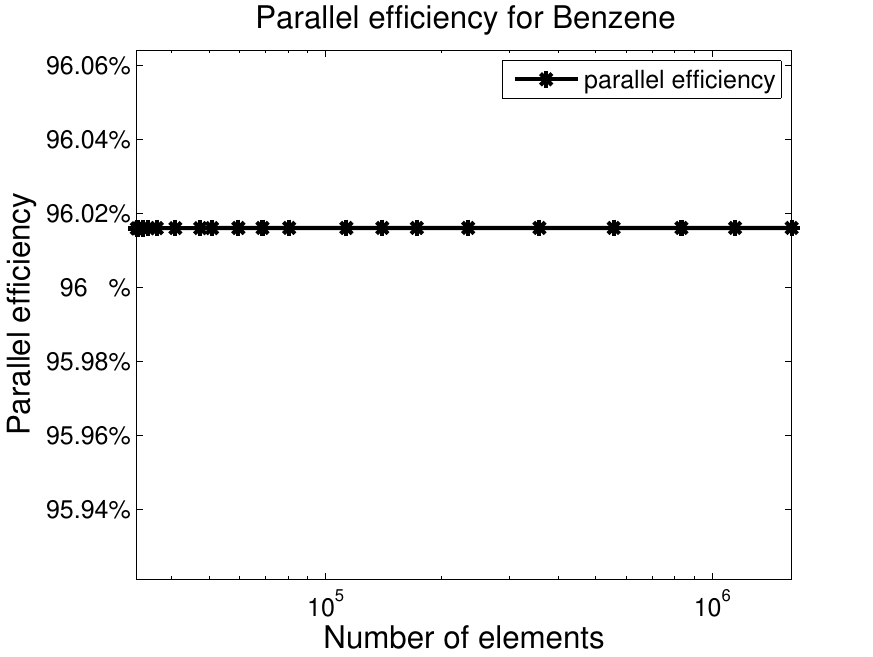}
\caption{Parallel efficiency of Algorithm \ref{multilevelAFEMim} for Lithium, Methane, Benzene and Ethanol.}
\label{ben-1}
\end{figure}

In the last subsection, we evaluate the parallel scalability of Algorithm \ref{multilevelAFEMim}. 
We calculate the trend of the parallel efficiency of Algorithm \ref{multilevelAFEMim} with the increase of the computing scale.
The corresponding results are presented in Figure \ref{ben-1}, which shows that quite a good parallel scalability can be obtained through solving
different wavefunctions by different computing nodes.

\section{Conclusion}
This paper presented an efficient multilevel correction type of adaptive algorithm for Hartree--Fock equation,
significantly improving computational efficiency while maintaining accuracy.
The method transforms the Hartree--Fock equation into tractable linearized boundary value problems and small-scale Hartree--Fock equations  in low-dimensional correction spaces,
avoiding direct computation of large-scale nonlinear eigenvalue problems and dense matrices.
Other key points include independent construction of correction spaces for parallelization, tensor-based preprocessing for computational reuse.
This approach enables efficient three-dimensional Hartree--Fock calculations and provides a solid foundation for accurate electronic structure simulations of complex systems.



\end{document}